\documentclass[12pt]{amsart}
\usepackage{amsfonts}
\usepackage{amsmath}
\usepackage{graphicx}
\usepackage{amssymb}
\newtheorem{theorem}{Theorem}[section]
\theoremstyle{plain}

\numberwithin{equation}{section}

\begin{document}
\title[Uniqueness of Viscosity Solutions]{Uniqueness of viscosity solutions of a geometric fully nonlinear parabolic equation}
\author{Jingyi CHEN}
\address{Department of Mathematics\\
University of British Columbia\\
Vancouver, B.C., V6T 1Z2\\
Canada}
\email{jychen@math.ubc.ca}
\author{Chao Pang}
\email{ottokk@math.ubc.ca}


\date{March 10, 2009}

\thanks{The first author is partially supported by an NSERC grant.}

\maketitle

{\small{{\bf{Abstract.}} 
We observe that the comparison result of Barles-Biton-Ley for viscosity solutions of 
a class of nonlinear parabolic equations can be applied to a geometric fully nonlinear parabolic equation 
which arises from the graphic solutions for the Lagrangian mean curvature flow. } }


\section{Introduction}

We consider the question of uniqueness for the following 
 fully nonlinear parabollic equation
\begin{equation}\label{eqn:1}
\frac{\partial u }{\partial t}=\sum\limits_{j=1}^{n}{\arctan \lambda_j}
\end{equation}
with initial condition 
$
u(x,0) = u_0(x), 
$
where $u$ is a function from ${\mathbb R}^n$ to ${\mathbb R}$ and $\lambda_j$'s are the eigenvalues of the Hessian $D^2u$. This equation arises naturally from geometry. In fact, when $u$ is a regular solution to (\ref{eqn:1}), it is known that the graph $(x,Du(x,t))$ evolves by the mean curvature flow and  it is a Lagrangian submanifold in ${\mathbb R}^n\times{\mathbb R}^n$ with the standard symplectic structure, for each $t$.  Recently, smooth longtime entire solution to (\ref{eqn:1}) has been constructed in \cite{CCH} assuming a certain bound on the Lipschitz norm of $Du_0$. 

Barles, Biton and Ley have obtained a very useful general  comparison result (Theorem 2.1 in \cite{BBL}) for the viscosity solutions for a class of fully nonlinear parabolic equations, as well as existence result (Theorem 3.1 in \cite{BBL}). In particular, they showed that  (\ref{eqn:1}) admits a unique longtime continuous viscosity solution for any continuous function $u_0$ in ${\mathbb R}$ when $n=1$. 

In this short note, we observe, via elementary methods, that the hypotheses in the general theorems in \cite{BBL} are valid for the geometric evolution equation (\ref{eqn:1}) in general dimensions. The result is the following

\begin{theorem}\label{main theorem}
Let $u$ and $v$ be an upper semicontinuous and a lower semicontinuous viscosity subsolution and supersolution  to (\ref{eqn:1}) in 
${\mathbb R}^n\times[0,T)$. If $u(x,0)\leq v(x,0)$ for all $x\in{\mathbb R}^n$, then $u\leq v$ in ${\mathbb R}^n\times[0,T)$. In particular, for any continuous function $u_0$ in ${\mathbb R}^n$, there is a unique 
continuous viscosity solution to (\ref{eqn:1}) in ${\mathbb R}^n\times[0,\infty)$. 
\end{theorem}

\section{Hypotheses (H1) and (H2)}

We now describe the assumptions in the comparison and existence results in \cite{BBL}. 

Let $S_n$ be the linear space of real $n\times n$ symmetric matrices. If $X\in S_n$, there exists an orthogonal matrix $P$ such 
that $X= P\Lambda P^T$ where $\Lambda$ is the diagonal matrix with diagonal entries consist of eigenvalues of $X$. Let  $\Lambda^+$ be the diagonal matrix obtained by replacing the negative eigenvalues in $\Lambda$ with 0's. Define $X^+=P\Lambda^+P^T$. 

Consider a continuous function $F$ from ${\mathbb R}^n\times[0,T]\times{\mathbb R}^n\times S_n$ to ${\mathbb R}$. The following assumptions of $F$ are necessary to apply the results in \cite{BBL}:

(H1) For any $R>0$, there exists a fuction $m_R : {\mathbb R}_+ \rightarrow {\mathbb R}_+ $ such that $m_R (0^+) = 0 $ and
$$
F(y,t,\eta(x-y),Y)-F(x,t,\eta(x-y),X)\le m_R(\eta \left|x-y\right|^2+\left|x-y\right|)
$$
for all $x,y \in \overline{\textit{B}}(0,R)$ and $t \in [0,T]$, whenever $X,Y \in \textsl{S}_n$ and $\eta > 0 $ satisfy

$$
-3\eta
\left(
\begin{array}{cc}
I & 0\\
0 & I
\end{array}
\right)
\le
\left(
\begin{array}{cc}
X & 0\\
0 & -Y
\end{array}
\right)
\leq
3\eta
\left(
\begin{array}{cc}
I & -I\\
-I & I
\end{array}
\right)
$$

(H2) There exists $0 < \alpha <1$ and constants $\textsl{K}_1>0$ and $\textsl{K}_2>0$ such that
$$
F(x,t,p,X)-F(x,t,q,Y)\le \textsl{K}_1\left|p-q\right|(1+\left|x\right|)+\textsl{K}_2\left(\mbox{tr}\,(Y-X)^+\right)^\alpha
$$
for every $(x,t,p,q,X,Y) \in {\mathbb R}^n\times [0,T] \times {\mathbb R}^n \times {\mathbb R}^n  \times \textsl{S}_n \times \textsl{S}_n$.

The operator $F$ is degenerate elliptic if (H2) holds. 

\begin{theorem}(Barles-Biton-Lay)\label{main}
Let u and $V$ be an upper semicontinuous viscosity subsolution and a lower semicontinuous viscosity supersolution respectively of
\begin{eqnarray*}
\frac{\partial u}
{\partial t} + F(x,t,Du,D^2u)&=&0 \,\,\,\,\,\, \,in \, \,\,\, {\mathbb R}^n \times [0,T)\\
u(\cdot,0)&=&u_0 \, \,\,\, in \,\, \,\,{\mathbb R}^n.
\end{eqnarray*}
Assume that (H1) and (H2) hold for $F$. Then

(1) If $u(\cdot,0)\le v(\cdot,0)$ in ${\mathbb R}^n$, then $u\le v$ in ${\mathbb R}^n \times  [0,T)$.

(2) If $u_0\in C({\mathbb R}^n)$ there is a unique continuous viscosity solution in ${\mathbb R}^n\times[0,\infty)$. 
\end{theorem}

We now present the proof of Theorem \ref{main theorem}. 

\begin{proof}

We define $F: \textsl{S}_n \rightarrow {\mathbb R}$ by 
\begin{equation}\label{eqn:2}
F(X)=-i\log \frac{\det(I+iX)}{{\det(I+X^2)}^{\frac{1}{2}}}
=-\frac{i}{2}\log \frac{\det(I+iX)}{\det(I-iX)}.
\end{equation}
That $F$ takes real values follows easily from 
\[
\overline{F(X)}=\frac{i}{2}\log \frac{\det(I-iX)}{\det(I+iX)}=F(X).
\]
Note that $F(D^2u)$, by diagonalizing $D^2u$ at a point, is equal to $\sum \arctan \lambda_j$. Therefore the flow (\ref{eqn:1}) can be written as 
 $$
 u_t + (-F(D^2u)) = 0.
 $$

Since $F(x,t,p,X)=F(X)$ is independent of $x$, the right hand side of the inequality for $F$ in (H1) must be zero, namely $m_R=0$. 
By multiplying an arbitrary vector $(\xi,\xi)\in {\mathbb R}^n\times{\mathbb R}^n$ and its transpose to the second matrix inequality in (H1), we see that  $X \leq Y$. Therefore, in order to establish (H1) it suffices to show: 

\vspace{0.15cm}

(H1') For any $X, Y \in \textsl{S}_n$, if $X \geq Y$  then $F(X) \geq F(Y)$.

\vspace{.15cm}

For any $X,Y\in S_n$ and $t\in[0,1]$, define 
$$
f_{XY}(t) = F(tX +(1-t)Y).
$$
We will show that $f_{\small{XY}}(t)$ is nondecreasing in $t\in[0,1]$ and then (H1') will follow as  $f_{XY}(0)=F(Y)$ and $f_{XY}(1)=F(X)$. 
Set  
$$A= I + i(tX+(1-t)Y)
$$ 
and 
$$B= I - i(tX+(1-t)Y).
$$ 
Then 
$$
f_{XY}(t) = -\frac{i}{2}(\log \det A - \log \det B).
$$
It follows that $AB=BA$ and 
$$
\left(A^{-1}+B^{-1}\right) \cdot \frac{AB}{2}= \frac{A+B}{2}=I. 
$$
Note that both $A$ and $B$ are invertible matrices for all $t\in[0,1]$. Hence, by using the formula $\partial _t \ln \det G = \mbox{tr} ( G ^{-1}\partial _t G)$ for $G(t)\in GL(n,{\mathbb R})$,  we have

\begin{eqnarray}\label{eqn:3}
f_{XY}'(t) & = &-\frac{i}{2} \mbox{tr}\left( A^{-1}\cdot\partial _t A-B^{-1}\cdot\partial _t B\right)\nonumber\\
& = & -\frac{i}{2} \mbox{tr}\left( (A^{-1}+B^{-1})\cdot i(X-Y) \right)\nonumber\\
& = & \mbox{tr} \left( (I+(tX+(1-t)Y)^2)^{-1} \cdot (X-Y)\right).
\end{eqnarray}
Since  $tX+(1-t)Y$ is real symmetric,  the matrix 
$$
C=I+(tX+(1-t)Y)^2
$$ 
is positive definite, hence so is  $C^{-1}$. There exists a matrix $Q\in GL(n,{\mathbb R})$ such that $C = QQ^T$. By the assumption $X \geq Y$, we have 
\begin{eqnarray*}
\mbox{tr} \left(C^{-1}(X-Y)\right) &=& \mbox{tr} \left( Q\cdot Q^T(X-Y) \right)\\
&=&\mbox{tr}\left(Q^T(X-Y)\cdot Q\right)\\
&\geq& 0
\end{eqnarray*}
since $Q^T(X-Y)Q$ is positive semidefinite. Therefore, we have shown that (H1) is valid for $F$ defined in (\ref{eqn:2}).





As $F(x,t,p,X)$ is independent of $p$, (H2) reads: there exist constants 
$K > 0$ and $0< \alpha <1$ such that $F(X)-F(Y)\leq K\left(\mbox{tr}(X-Y)^+\right)^\alpha$ for all $X,Y\in S_n$.

For any $X, Y\in S_n$,  integrating (\ref{eqn:3})  leads to  
\begin{equation}\label{eqn:4}
F(X)-F(Y)=\int_0^1 \mbox{tr}\left(C^{-1}X\right)\,dt. 
\end{equation}

For $X-Y\in S_n$ there exists an orthogonal matrix $P$ such that $X-Y=P\Lambda P^T$ where the diagnal matrix 
$\Lambda$ has diagonal entries $\lambda_1, ... , \lambda_n$. Let $\lambda^+_j=\max\{\lambda_j,0\}$. Since $0<C^{-1} \leq I$, we have $0<P^TC^{-1}P\leq I$. If $c_{jj}$ denote the diagonal entries of $P^TC^{-1}P$ for $j=1, ..., n$, 
then $c_{jj} = \langle P^TC^{-1}P e_j,e_j\rangle$ where $\{e_1, ... , e_n\}$ is the standard basis for ${\mathbb R}^n$ and $\langle\cdot,\cdot\rangle$ is the Euclidean inner product. It follows that $0<c_{jj}\leq 1$ for $j=1, ... , n$. 
 Then 
\begin{eqnarray*}
\mbox{tr}\left(C^{-1}(X-Y)\right) &=& \mbox{tr}\left( P^T C^{-1}P \cdot P^T(X-Y)P \right)\\
&=&\mbox{tr}\left( P^TC^{-1}P\cdot\Lambda\right)\\
&=&\sum c_{jj}\lambda_j\\
&\leq& \sum \lambda^+_j\\
&=&\mbox{tr} (X-Y)^+.
\end{eqnarray*}
Substituting the above inequality into (\ref{eqn:4}) implies: for any $X,Y\in S_n$ we have
$$
F(X)-F(Y)\leq \mbox{tr} \left(X-Y\right)^+.
$$

Because $\arctan x$ is in $(-\pi/2,\pi/2)$, we have $F(X)-F(Y)<n\pi$. For any constant $\alpha$ with $0<\alpha<1$, 
if $\mbox{tr}\left(X-Y\right)^+ \leq 1$ then 
$$
F(X)-F(Y)\leq \mbox{tr}\left(X-Y\right)^+\leq n\pi\, \mbox{tr}\left[\left(X-Y\right)^+\right]^\alpha
$$ 
and if $\mbox{tr}\left(X-Y\right)^+>1$ then 
$$
F(X)-F(Y)\leq n\pi \leq n\pi \, \mbox{tr}\left[\left(X-Y\right)^+\right]^\alpha.
$$ 
Therefore, (H2) holds for $K_2=n\pi$ and any constants $K_1>0$ and  $\alpha$ with $0<\alpha<1$. 

Now Theorem \ref{main theorem} follows immediately from Theorem \ref{main}. \end{proof}

We remark on that  (H1'), for the operator $F(X)=\sum\arctan\lambda_j(X)$, also follows from the basic fact (cf. p.182 in \cite{HJ}):  Suppose that $X,Y\in S_n$ and the eigenvalues $\lambda_j$'s of $X$ and $\mu_j$'s of $Y$ are in descending order $\lambda_1\geq\lambda_2\geq ... \geq\lambda_n$ and $\mu_1\geq\mu_2\geq ... \geq\mu_n$. If $X\geq Y$, then $\lambda_j\geq \mu_j$ for $j=1, ... , n$. 

We also mention the uniqueness of viscosity solutions of the Cauchy-Dirichlet problem for (\ref{eqn:1}). 
Note that the operator $F(X)=\sum\arctan\lambda_j(X)$ satisfies (H1') which is exactly the fundamental monotonicity condition (0.1) for $-F$ in \cite{CIL}, therefore $-F$ is proper in the sense of \cite{CIL} (cf. p.2 in \cite{CIL}). As (H1) holds, Theorem 8.2 in \cite{CIL} is valid for (\ref{eqn:1}): 

\begin{theorem}
The continuous viscosity solution to the following Cauchy-Dirichlet problem is unique:
\begin{eqnarray*}
&&  u_t  = \sum^n_{j=1}\arctan\lambda_j, \,\,\,\,\,\mbox{in $ (0,T) \times \Omega$} \\
&& u(t,x) = 0, \,\,\,\,\,\,\,\mbox{ for $0 \leq t<T$ and $x \in \partial \Omega$} \\
&& u(0,x) = \psi (x), \,\,\,\,\mbox{for $x \in \overline{\Omega}$}
\end{eqnarray*}
where $\lambda_j$'s are the eigenvalues of $D^2u$, $\Omega \subset {\mathbb R}^n$ is open and bounded and $T > 0$ and $\psi \in C(\overline{\Omega})$. If $u$ is an upper semicontinuous viscosity solution and $v$ is a lower semicontinuous viscosity solution of the Cauchy-Dirichelt problem, then $u\leq v$ on $[0,T)\times\Omega$. 
\end{theorem}
Note that the initial-boundary conditions for the subsolution and supersolution are: $u(x,t)\leq 0 \leq v(x,t)$ for $t\in[0,T)$ and $x\in\partial\Omega$ and $u(x,0)\leq \psi(x)\leq v(v,0)$ for $x\in \overline{\Omega}$.


\begin{thebibliography}{99}
\bibitem{BBL}G. Barles, S. Biton, O. Ley; {\it Uniqueness for parabolic equations without growth condition and applications to the mean curvature flow in $\textbf{R}^2$,} J. Differential. Equations  {\bf 187}  (2003), 456-472.
\bibitem{CCH} A. Chau, J. Chen, W. He; {\it Lagrangian mean curvature flow for entire Lipschitz graphs}, preprint. 
\bibitem{CIL}M.G. Crandall, H. Ishii, P.-L. Lions; {\it User's guide to viscosity solutios of second ordr partial differential equations,} Bull. Amer. Math. Soc  {\bf 27}  (1992),  1-67.
\bibitem{HJ} R.A. Horn, C.R. Johnson, Matrix Analysis, Cambridge University Press, 1985. 
\end{thebibliography}
\end{document}